\documentclass{amsart}
\usepackage{amssymb}
\usepackage{enumerate}
\usepackage{url}

\newcommand\N{\mathbb N}

\newcommand\R{\mathbb R}

\newcommand\ph\varphi
\newcommand\ps\psi
\newcommand\ep\varepsilon
\newcommand\rh\varrho
\newcommand\al\alpha
\newcommand\be\beta
\newcommand\ga\gamma
\newcommand\om\omega
\newcommand\ta\tau
\renewcommand\th\theta
\newcommand\de\delta
\newcommand\ze\zeta
\newcommand\ch\chi
\newcommand\et\eta
\newcommand\io\iota
\newcommand\la\lambda
\newcommand\si\sigma

\newcommand\Ga\Gamma
\newcommand\De\Delta
\newcommand\Th\Theta
\newcommand\La\Lambda
\newcommand\Si\Sigma
\newcommand\Ph\Phi
\newcommand\Ps\Psi
\newcommand\Om\Omega

\newtheorem{theorem}{Theorem}
\newtheorem{lemma}[theorem]{Lemma}

\newtheorem{corollary}[theorem]{Corollary}
\theoremstyle{definition}

\theoremstyle{remark}

\newcommand\x{\bar X}

\begin{document}
\title[Positive polynomials with structured sparsity]
{A note on the representation of positive polynomials with structured sparsity}

\author{David Grimm}
\address{Universit\"at Konstanz, Fachbereich Mathematik und Statistik, 78457 Konstanz, Germany} 
\email{david.grimm@uni-konstanz.de}

\author{Tim Netzer}
\email{tim.netzer@uni-konstanz.de}

\author{Markus Schweighofer}
\email{markus.schweighofer@uni-konstanz.de}

\date{\today}

\begin{abstract}
We consider real polynomials in finitely many variables. Let the variables
consist of finitely many blocks that are allowed to overlap in a certain way.
Let the solution set of a finite system of polynomial inequalities be
given where each inequality involves only variables of one block. We
investigate polynomials that are positive on such a set and sparse in
the sense that each monomial involves only variables of one block.
In particular, we derive a short and direct proof for Lasserre's theorem
of the existence of sums of squares certificates respecting the block
structure. The motivation for the results can be found in the literature
and stems from numerical methods using semidefinite programming to simulate or
control discrete-time behaviour of systems.
\end{abstract}

\maketitle

\section{Introduction}
Let $\R[\x]:=\R[X_1,\ldots,X_n]$ be the ring of real
polynomials in $n$ variables $\x:=(X_1,\dots,X_n)$. A subset
$M\subseteq\R[\x]$ is called a \emph{quadratic module} of $\R[\x]$ if it
it contains $1$ and is closed under addition and multiplication with squares,
i.e., $1\in M$, $M+M\subseteq M$ and $\R[\x]^2M\subseteq M$. A quadratic module
$M\subseteq\R[\x]$ is called \emph{archimedean} if
$N-\sum_{i=1}^nX_i^2\in M$ for some $N\in\N$.

Let $g_1,\ldots, g_m\in\R[X]$ be given and set
$g_0:=1\in\R[\x]$. These polynomials define a
\emph{basic closed semialgebraic set}
$$S:=
\{x\in\R^{n}\mid g_{1}(x)\geq 0,\ldots, g_m(x)\geq 0\}\subseteq\R^n$$
and generate in $\R[\x]$ the quadratic module
$$M:=
\left\{\sum_{i=0}^m\si_ig_i\mid\si_0,\dots,\si_m\in\sum\R[\x]^2\right\}
\subseteq\R[\x]$$
where $\sum\R[\x]^2$ denotes the set of sums of squares in the ring $\R[\x]$.
Obviously, every polynomial from $M$ is nonnegative on $S$.
In particular, if $M$ is archimedean, then $S$ is compact.
In 1991, Schm\"udgen \cite{s} proved a very surprising partial converse to
these facts (for example, if
$MM\subseteq M$, then $M$ is archimedean if and only if $S$ is compact).
See \cite[Section 1]{sch} and \cite{pd} for this and a discussion when $M$ is archimedean.
Shortly after this groundbreaking work of Schm\"udgen, Putinar \cite{put}
proved the following theorem (note that $f>0$ on $S$ should mean that $f$ is
\emph{strictly} positive on $S$):

\begin{theorem}[Putinar]\label{putinar}
Suppose $M$ is archimedean. Then for every $f\in\R[\x]$,
$$\text{$f>0$ on $S$}\implies f\in M.$$
\end{theorem}

A very special case of this is the following theorem of Cassier
\cite[Th\'eor\`eme 4]{cas} (which can even be derived from earlier results
like \cite[Th\'eor\`eme 12]{kri}).

\begin{corollary}[Cassier]\label{cassier}
For all $R>0$ and $f\in\R[\x]$, if $f>0$ on the closed ball centered at the
origin of radius $R$, then there are $\si,\ta\in\sum\R[\x]^2$ such that
$f=\si+\ta(R^2-\sum_{i=1}^n X_i^2).$
\end{corollary}

In 2001, Lasserre \cite{l1} recognized that Putinar's result can be used
for numerical computation of the minimum of a polynomial on a non-empty
\emph{compact} basic closed semialgebraic $S$ by semidefinite programming (a
well-known generalization of linear programming).
The idea is to maximize $\la$ such that $f-\la$ lies in a certain
finite-dimensional subset of $M$ that can be expressed in a semidefinite
program (SDP for short) whose size depends on the size of the chosen subset
of $M$.
Since $M$ can be exhausted by such subsets, the results of Schm\"udgen and
Putinar say that the accuracy of this method is arbitrarily good for large
SDPs \cite{l1,sch}.

In many problems, the polynomials $f$ and $g_i$ are sparse.
Waki, Kim, Kojima and Muramatsu tried to take advantage of this sparsity and
implemented corresponding SDP relaxations that turned out to be very
efficient in practice \cite{w}. They correspond to a certain subset of the
quadratic module $M$ reflecting the sparsity pattern of $f$ and the $g_i$.
But a theoretical result on the accuracy of such ``sparse'' SDP
relaxations was only given recently by Lasserre \cite{l2}. Lasserre's proof
uses SDP duality, compactness arguments,
theorems about the construction of measures with given marginals and
Putinar's solution to the moment problem (a dual result to Theorem
\ref{putinar}). Moreover, Lasserre assumes that $S$ has non-empty interior.
For general compact $S$, the result follows from an additional argument
given by Kojima and Muramatsu \cite{km}.

In this note, we give a short direct proof of Lasserre's result taking only
Corollary \ref{cassier} for granted. 

\section{The Theorem}

For $I\subseteq \{1,\ldots,n\}$, we denote by $X_{I}$ the set of variables
$\{ X_{i}\mid i\in I \}$ and by $\R[X_{I}]$ the polynomial ring in these
variables. By $\R^{I}$ we denote the subspace of $\R^{n}$ which corresponds
to the variables in $X_{I}$. In particular, an element of $\R[X_{I}]$ can be
seen as a function on $\R^{I}$.

Following \cite{l2}, we suppose from now on that we have non-empty
$I_{1},\ldots,I_{r}\subseteq \{1,\ldots,n\}$
satisfying the following \textit{running intersection property:}
\begin{equation}\tag{RIP}\label{rip}
\text{For all $i=1,\dots,r$, there is $k<i$ such that
$I_{i}\cap \bigcup_{j<i}I_{j} \subseteq I_k$}.
\end{equation}
Note that \eqref{rip} is empty in the case $r\le 2$.

\begin{lemma}\label{box}
Let $K\subseteq\R$ be compact. Suppose $f=f_{1}+\ldots + f_{r}$ with
$f_{j}\in\R[X_{I_{j}}]$ and $f>0$ on $K^n$.
Then $$f=h_{1}+\ldots +h_{r}$$
for some $h_{j}\in\R[X_{I_{j}}]$ with $h_j>0$ on
$K^{I_{j}}$.
\end{lemma}

\begin{proof}
The proof is by induction on $r$. First suppose $r=2$. We may of course assume $K\neq\emptyset$.
There is some $\ep>0$ such that $f=f_1+f_2\ge \ep$ on $K^n$. 
Now consider the function $h:K^{I_{1}\cap I_{2}}\to\R$ defined by
$$h(y)=\min\{f_{1}(x,y)\mid x\in K^{I_1\setminus I_2}\}-\frac{\ep}{2}.$$ Note that $K^{\emptyset}$ is a singleton and thus $h+\frac{\ep}{2}$ is the minimum of $f_{1}$ on $K^{I_{1}}$ in the case $I_1\cap I_2=\emptyset$. 
To show that $h$ is continuous, consider
$y,y'\in K^{I_1\cap I_2}$. Choose
$x,x'\in K^{I_1\setminus I_2}$ minimizing $f(x,y)$ and $f(x',y')$,
respectively.
Then
\begin{align*}
|h(y)-h(y')|&=|f_1(x,y)-f_1(x',y')|\\
&\le\max\{|f_1(x,y)-f_1(x,y')|,|f_1(x',y)-f_1(x',y')|\}
\end{align*}
shows that $h$ is uniformly continuous because $f_1$ is uniformly
continuous on the compact set $K^{I_1}$. We now claim that
$$f_1-h\ge\frac{\ep}{2}\text{\ on\ }K^{I_1}\qquad\text{and}\qquad
  f_2+h\ge\frac{\ep}{2}\text{\ on\ }K^{I_2}.$$
The first claim is clear by the definition of $h$. To show the second,
suppose $(y,z)\in K^{I_{2}}=K^{I_{1}\cap I_{2}}\times
K^{I_{2}\setminus I_{1}}$, choose
$x\in K^{I_{1}\setminus I_{2}}$ with
$h(y)=f_1(x,y)$ and observe
$$f_2(y,z)+h(y)=f_2(y,z)+f_1(x,y)-\frac{\ep}{2}=f(x,y,z)-\frac{\ep}{2}\ge\frac{\ep}{2}.$$
Now approximate $h$ by a polynomial $p\in\R[X_{I_{1}\cap I_{2}}]$ such that
$|h-p|\le\frac\ep 4$ on $K^{I_{1}\cap I_{2}}$.
Then 
$$h_1:=f_1-p > 0 \text{\ on\ }K^{I_1}\qquad\text{and}\qquad
  h_2:=f_2+p > 0 \text{\ on\ }K^{I_2}$$
and $f=h_{1}+h_{2}$. This completes the proof for $r=2$.

For the induction step, suppose $r\ge 3$. Setting
$$\tilde{f}:= f_{1} +\ldots +f_{r-1}\in \R[X_{\bigcup_{j<r}I_{j}}],$$
we have $f=\tilde{f}+f_{r}$. Now the proof for the case $r=2$ showed that
there is some polynomial $p$, using only the variables indexed in
$I_{r}\cap \bigcup_{j<r}I_{j}$, such that $\tilde{f}-p$ and $f_{r}+p$ are
positive on appropriate cartesian powers of $K$. By \eqref{rip},
$p\in\R[X_{I_{k}}]$ for some $k<r$. Hence $\tilde{f}-p$ is
a sum of polynomials from the rings
$\R[X_{I_{j}}], j=1,\ldots, r-1$. The induction hypothesis therefore applies
to $\tilde{f}-p$ and yields the result.
\end{proof}

Again following \cite{l2}, we now suppose that for $j=1,\ldots, r$, we have
polynomials $g_1^{(j)},\ldots,g_{l_{j}}^{(j)}\in\R[X_{I_{j}}]$ which
define sets
$$S_{j}:=\{x\in\R^{I_{j}}\mid
g_{1}^{(j)}(x)\ge 0,\dots,g_{l_{j}}^{(j)}(x)\ge 0\}\subseteq \R^{I_{j}}.$$ We do not require these sets to be compact at the moment.
Let $S\subseteq\R^n$ be the basic closed semialgebraic set defined by all
these polynomials $g_{i}^{(j)}$ in $\R^n$.

\begin{theorem}\label{thm}
Suppose $f=f_1+\ldots+f_r$ with
$f_j\in\R[X_{I_j}]$ and $f>0$ on $S$. Then for any bounded set
$C\subseteq\R^n$, there
are $0<\la\le 1$, $k\in\N$ and polynomials $h_{j}\in\R[X_{I_{j}}]$
with $h_j>0$ on $C$ such that
$$f=\sum_{j=1}^{r}\sum_{i=1}^{l_{j}} 
(1-\lambda g_{i}^{(j)})^{2k}g_{i}^{(j)} + \sum_{j=1}^{r}h_{j}.$$
\end{theorem}

\begin{proof}
Choose a compact set $K\subset\R$ such that $C\subseteq K^n$.
Choose $0<\la\le 1$ such that  $\la g_{i}^{(j)}\le 1$ on $K^{n}$ for all
$i,j$. 
Set $$f_k := f-\sum_{j=1}^{r}\sum_{i=1}^{l_{j}}\left(1-\la g_{i}^{(j)}\right)^{2k}g_{i}^{(j)} \qquad (k \in \N).$$
We have $f_k \leq f_{k+1}$ on $K^n$ for $k\in \N$ and one checks that for all $x \in K^n$ there exists $k \in \N$ such that $f_k(x)>0$. By compactness of $K^n$, we find therefore $k \in \N$ with $f_k>0$ on $K^n$.   
Now apply Lemma
$\ref{box}$ to $f_k$.
\end{proof}

For $S$ with non-empty interior, the following is Lasserre's result
\cite[Corollary 3.9]{l2}. It has been extended to general compact $S$
by Kojima and Muramatsu in \cite[Theorem 1]{km}.

\begin{corollary}[Lasserre, Kojima, Muramatsu]\label{cor}
Let the quadratic modules $M_j$
generated by $g_1^{(i)},\dots,g_{l_j}^{(j)}$ in $\R[X_{I_{j}}]$
be archimedean. If $f\in\R[X_{I_1}]+\dots+\R[X_{I_r}]$ and
$f>0$ on $S$, then $f\in M_1+\dots+M_r$.
\end{corollary}

\begin{proof}
Choose $R>0$ big enough such that $R^2-\sum_{i\in I_j}X_i^2\in M_j$ for
all $j$. Define $C\subseteq\R^n$ to be the closed ball with radius $R$
around the origin. Now apply Theorem \ref{thm} and note that $h_j\in M_j$
by Corollary \ref{cassier}.
\end{proof}

\end{document}